\documentclass{amsart}
\usepackage{psfig,latexsym,amssymb}
\def\proof{\medskip\noindent{\sc Proof. }}
\def\proofo{\medskip\noindent{\sc Proof of Theorem 2.1. }}
\def\EOP{\hfill$\Box$}

\newtheorem{thm}{Theorem}[section]
\newtheorem{lem}{Lemma}[section]

\newtheorem{defn}{Definition}[section]

\newtheorem{conj}{Conjecture}[section]

\begin{document}

\title{A Hodge decomposition for the complex of injective words}
\author{Phil Hanlon}
\author{Patricia Hersh}
\address{Department of Mathematics,
University of Michigan, Ann Arbor, Michigan 48109-1109}
\curraddr{
        Department of Mathematics, 
        University of Michigan,
        525 East University Ave.
        Ann Arbor, Michigan 48109-1109}
\thanks{The second author was supported by an NSF postdoctoral 
research fellowship.}
\begin{abstract}
Reiner and Webb compute the $S_n$-module structure for the complex of
injective words in [RW].
This paper refines their formula by providing a Hodge type decomposition.
Along the way, this paper proves
that the simplicial boundary map interacts in a nice fashion with the 
Eulerian idempotents.  

The Laplacian acting on the top chain group in the complex of injective
words is also shown to equal the signed random to random shuffle operator.
Uyemura-Reyes conjectures in
[Uy] that the (unsigned) random to random shuffle operator has integral
spectrum.  We prove that this conjecture would imply that
the Laplacian on (each chain group in) the complex
of injective words has integral spectrum.
\end{abstract}

\maketitle

\section{Introduction}\label{introsec}
Let $V = \langle v_1, \ldots ,v_n \rangle $ be 
an n-dimensional Euclidean space.
For each $r$, let $\Gamma_r = V^{\otimes r}$ and let $\partial
_r : 
\Gamma_r \rightarrow \Gamma_{r-1}$ be the map given by:
$$\partial
_r(a_1 \otimes \cdots \otimes a_r) = \sum_{j=1}^r 
(-1)^{j-1} (a_1 \otimes \cdots \otimes a_{j-1} \otimes a_{j+1} \otimes
\cdots \otimes a_r).$$
It is well-known that the $\partial
_r$ are boundary maps, i.e., that $
\partial
_r \cdot \partial
_{r+1} = 0$.  
 
Let $M_r$ be the multilinear part of $\Gamma_r$.  So, $M_r = 0$ if 
$r > n$, and for $r \le n$,
 
$$M_r 
= \langle v_{i_1} \otimes v_{i_2} \otimes \cdots \otimes v_{i_r} : 
i_1, \ldots ,i_r \hspace{.05in}{\rm are} \hspace{.05in} {\rm distinct } 
\rangle .$$ 
Note that $\dim (M_r) = n(n-1) \cdots (n-r+1) = {{n!} \over {(n-r)!}}$.  
Also, it is clear that $\partial
_r (M_r) \subset M_{r-1}$ and so
$$0 \rightarrow M_n \rightarrow M_{n-1} \rightarrow \cdots \rightarrow M_0
\rightarrow 0$$
is a subcomplex of $(\Gamma_*, \partial
_*)$.  This paper will concern the homology of this
subcomplex.
 
The complex $(M_*,\partial
_*)$ appears in earlier work on the subword order of
injective words on the alphabet $\{1,2,\ldots,n\}$.  This poset is the face
poset of a regular CW complex $K_n$ whose homology
agrees with the homology
of $(M_*,\partial
_*)$.  In [F], Farmer proves that $K_n$ is homotopy equivalent
to a wedge of $(n-1)$-spheres thus showing that the homology of $(M_*,
\partial
_*)$ vanishes except at top degree.  Bjorner and Wachs [BW] prove a
stronger result -- that the lexicographic order on permutations induces a
recursive coatom ordering on the poset of injective words.  This in turn 
gives a dual CL-shelling of $K_n$.
 
Reiner and Webb [RW] study $(M_*,\partial
_*)$ as a subcomplex of $(\Gamma_*,
\partial
_*)$.  The natural action of $S_n$ on $\{v_1, \ldots ,v_n \}$ extends
to an action of $S_n$ on $(\Gamma_*,\partial
_*)$ which preserves $(M_*,\partial
_*)$.
Reiner and Webb compute the homology of $(M_*,\partial
_*)$ as an $S_n$-module.
\vskip .3 cm
\noindent
\begin{thm}[Reiner-Webb]  As an $S_n$-module, the top homology of 
$(M_*,\partial_*)$ is 
$$\bigoplus_{k=0}^n (-1)^{n-k} ind_{S_{n-k}}^{S_n} (\varepsilon _{n-k})
= \bigoplus_{k=0}^n (-1)^{n-k} ind_{S_{n-k}\times S_k}^{S_n} 
(\varepsilon_{n-k} \otimes Reg_k)
$$
where $\varepsilon _{n-k}$ denotes the trivial representation of $S_{n-k}$.
Furthermore, the multiplicity of an irreducible $S^{\lambda}$ of 
$S_n$ in the top 
homology is equal to the number of standard Young tableaux of shape
$\lambda$ which have their smallest descent even.
\end{thm}
 
In this paper, we will do two
things.  First, we will show that there is a
natural {\it Hodge decomposition} of the homology of $(M_*,\partial_*)$.  
This decomposition will split $H_n(M)$ into $n$ components
$$H_n(M) = \bigoplus_{j = 1}^n H_n^{(j)}(M).$$
We will show that the dimension of $H_n^{(j)}(M)
$ is equal to the number of 
derangements with exactly $j$ cycles.  More specifically, we will 
show that each $H_n^{(j)}(M)
$ is invariant under the action of $S_n$ and prove that
$H_n^{(j)}(M)
$ is a sum of linear characters induced from centralizers of 
permutations with exactly $j$ cycles.
 
Second, we will study the Laplacian $\Lambda_*$ associated to the complex
$(M_*,\partial_*)$.  
We will show that $\Lambda_n$ is closely connected to the 
transition matrix for random to random shuffling.  Random to random shuffling
has been studied by Uyemura-Reyes in [Uy].  In [Uy], the author 
makes a conjecture about the spectrum of the transition matrix for random
to random shuffling which together with our results imply the conjecture 
that the spectrum of $\Lambda_n$ is integral.  We go on to compute 
$\Lambda_r$, for $r<n$ in terms of $\Lambda_n$.  This computation shows
that $\Lambda_r$ is positive definite for $r<n$ thus giving another proof
that $H_r(M) = 0$ for $0 \le r < n$.  This computation also shows that 
the spectrum of $\Lambda_r$ is integral if the spectrum of $\Lambda_n$ 
is integral.  Thus, if Uyemura-Reyes' conjecture 
on the spectrum of random to random
shuffling is correct, then the spectra of all $\Lambda_r$ are integral. 
 
\section{A Hodge type decomposition of $H_*(M)$}

We begin by recalling the definition of the Eulerian idempotents 
$e_r^{(j)}$ in 
{\bf C}$S_r$.  For each $r$ and $k$, let $S(r;k)$ denote the set of
permutations in $S_r$ with exactly $k-1$ descents.  Following Loday
[L], define elements $l_r^{(k)}$ and $\lambda_r^{(k)}$ in {\bf C}$S_r$
according to the following formulae:
$$ (2.1)\ \ \ 
l_r^{(k)} = (-1)^{k-1} \sum_{\sigma \in S(r;k)} sgn (\sigma) \sigma$$
$$ (2.2)\ \ \ 
\lambda_r^{(k)} = \sum_{i=0}^{k-1} (-1)^i {{n+i} \choose i} l_r^{k-i}$$
$$ (2.3)\ \ \ 
(-1)^{k-1} \lambda_r^{(k)} = \sum_{j=1}^n k^j e_r^{(j)}.$$
It is worth noting that the first two equations define the 
$\lambda_r^{(k)}$ explicitly.  The third equation then determines 
the $e_r^{(j)}$ in terms of the $\lambda_r^{(k)}$ because the transition
matrix $(k^j)_{k,j}$ is a Vandermonde and hence invertible.
 
There is a significant literature on the Eulerian idempotents and their 
many remarkable properties.  We will need two of these properties.  The
first is the well-known fact that $e_n^{(1)}, e_n^{(2)}, \ldots ,
e_n^{(n)}$ form a set of pairwise orthogonal idempotents in {\bf C}$S_n$.
In other words, $e_n^{(j)}\cdot e_n^{(\ell)} = 0$ if $j \not= \ell$ and 
$e_n^{(1)} + e_n^{(2)} + \cdots + e_n^{(n)} = id$.  This
implies that if $X$ is any $S_n$-module, then
$$X = \bigoplus_j 
e_n^{(j)}\cdot X.$$
The second fact we will need describes the relationship between the 
$e_n^{(j)}$ and the boundary map $\partial
$.  To state this result, it will be
helpful to write permutations in one-line notation.  Let $i \in \{1,2,
\ldots ,n \}$ and let $S_{n\backslash \{i\}}$ denote permutations of 
$\{1,2,\ldots ,n \} \backslash \{i\}$.  There is a natural identification 
of $S_{n \backslash \{i\}}$ with $S_{n-1}$ which comes about by changing 
each occurrence of $i+j$ to $i+j-1$.  Via this identification, we can think 
of $e_{n-1}^{(\ell)}$ as sitting inside of the group algebra of 
$S_{n \backslash \{i\}}$.  
 
As in [RW], we will think of $\partial
$ as acting on linear combinations of 
injective words on the alphabet $\{1,2,\ldots ,n\}$.  For $i \in 
\{1,2, \ldots ,n\}$, let $\partial
 [i]$ denote the part of $\partial
$ which removes
the number $i$.  The next theorem presents a surprisingly elegant outcome
to the computation of $\partial
 [i] e_n^{(k)}$.
 
Before stating and proving the result, we will give an example to be
sure that the notation is clear.  Let $n=3$.  In the example that follows, 
we will use $a,b,c$ in place of $1,2,3$ so as to avoid confusion with
coefficients.  The Eulerian idempotents are given by:
$$e_3^{(1)} = {1 \over 6}( 2 \cdot abc + bac + acb - 
  bca - cab - 2 \cdot cba) $$
$$e_3^{(2)} = {1\over 2} ( abc + cba) $$
$$e_3^{(3)} = {1 \over 6}( abc - bac - acb + bca + cab - cba ) $$

For the purposes of this example, we will apply $\partial [b]$
to each of these.  Doing so, we get
 
$$\partial
[b]e_3^{(1)} = {1 \over 6}( -2ac + ac + ac -ca -ca + 2ca) = 0$$
$$\partial
[b]e_3^{(2)} = {1 \over 2}( -ac -ca) = -e_2^{(1)}$$
$$\partial
[b]e_3^{(3)} = {1 \over 6}( -ac-ac-ac+ca+ca+ca) = -e_2^{(2)}.$$
\vskip .1 cm
\noindent

Part 2 of Theorem 2.1 is needed later to show that the complex of 
injective words has a Hodge decomposition.

\begin{thm}
Fix $n$ and $i \in \{1,2, \ldots ,n \}$.  Then,
\begin{enumerate}
\item
$\partial
[i]\lambda_n^{(k)} = (-1)^{i-1}k \lambda_{n-1}^{(k)}$
\item
$\partial
[i]e_n^{(k)} = (-1)^{i-1} e_{n-1}^{(k-1)}.$
\end{enumerate}
\end{thm}

Before proving this theorem, let us verify two lemmas to be used
in its proof.

\begin{lem}\label{descent-increase}
Among the $n$ places the letter $i$ could be inserted into a permutation
$\sigma\in S_{[n]\setminus \{ i\} }$ which has $j$ descents,
$j+1$ choices yield 
permutations with $j$ descents, while the other $n-j-1$
choices all yield permutations with $j+1$ descents.  
\end{lem}

\proof
First we consider the case $i=n$, then use a graph for a permutation to
generalize to all $i$.
Notice that inserting $n$ between two letters descending
letters preserves the number of descents, while inserting $n$ 
between two ascending letters increases the number of descents.  Thus, there
are $n-2-j$ ways to increase the number of descents by one by inserting
$n$ between two ascending letters.  In addition, inserting $n$ before the 
first letter gives one more way to increase the number of descents by one.

For $i\ne n$, the analysis will also need to consider whether or not $i$ is
intermediate in value to the pair of consecutive labels where it is to be
inserted.  To this end, we define the graph of a permutation as follows.

\begin{defn}
For each $\pi\in S_{[n]\setminus \{ i\} }$, define the related function
$\pi ' : [0,n] \rightarrow [0,n+1] $
by $\pi '(j) = \pi (j) $ for $1\le j < i$ and $\pi '(j) = \pi (j+1)$ for 
$n-1\ge j>i$.
In addition, let $\pi '(0)=n+1$ and $\pi '(n) = 0$.
Then the {\it graph} of $\pi $ is
obtained by plotting the points $(j,\pi '(j))$ for each $j\in [0,n] $, 
and for each $j \in [0,n-1]$ 
connecting the point $(j,\pi '(j))$ to $(j+1,\pi '(j+1))$ 
by a straight line segment.  
\end{defn} 

This graph has negative slope at each descent and positive slope at each
ascent.
Furthermore, it crosses the line $y=i$ with
negative slope one more time than 
it does with positive slope, because it represents a 
continuous function which begins above the line $y=i$ and ends below the 
line $y=i$ (and which has nonzero slope everywhere it touches the line
$y=i$).  

We claim that the number of places to insert $i$ which will increase
the number of descents by one is equal to sum of the number of 
ascents that do not cross the line together with the number of descents which
do cross the line.  This is clear except at the endpoints.  
By letting $\pi '(0)=n+1$ and
$\pi '(n+1) = 0$, we created descents at the initial and final positions in
$\pi $, which are only counted above when the graph crosses the line $y=i$
at these points, namely when $\pi (1) < i$ and when $\pi (n) > i$,
respectively.  These are exactly the situations where inserting $i$ at the 
initial or final positions will indeed increase the number of descents by
one.

We already observed that the number of descents crossing
the line is one more than the number of ascents crossing the line.
Thus, the total number of ways to increase the number of descents by one 
is one more than the total number of ascents, i.e. it is $n-j-1$, as desired.
A similar argument shows that all of the remaining $j+1$ options will 
preserve the number of descents.
\EOP

\begin{lem}\label{sign-change}
If $\tau \in S_n$ is obtained from $\sigma \in S_{n\setminus \{ i\} }$ 
by inserting $i$ after the $(d-1)$-st letter of $\sigma $, then
$sgn (\tau ) = (-1)^{i-d} sgn (\sigma )$.  
\end{lem}

\proof
If $d=1$, so $i=\tau_1$, then
inserting $i$ created $i-1$ new descents, because the values $1,\dots ,i-1$
all appear later than the value $i$.  Thus, 
$sgn (\tau ) = (-1)^{i-1} sgn (\sigma )$ in this case, as desired.  
Now we proceed by induction on $d$.  Moving the letter $i$ from position 
$r$ to position $r+1$ in $\tau $ by an adjacent transposition 
reverses the sign of $\tau $.  
Likewise, increasing $d$ from $r$ to $r+1$ reverses the 
sign of $(-1)^{i-d} sgn(\sigma )$ 
from $(-1)^{i-r} sgn (\sigma )$ to $(-1)^{i-(r+1)} sgn (\sigma )$, so 
$sgn (\tau )$ continues to agree with $(-1)^{i-1} sgn (\sigma )$ as 
$d$ increases.
\EOP

\proofo
We first prove identity (1).  Note that 
$$\lambda_n^{(k)} = \sum_{i=0}^{k-1} (-1)^i {{n+i} \choose i} 
(-1)^{k-i-1} \sum_
{\sigma} sgn(\sigma) \sigma$$where the sum is over
$\sigma$ with $k-i-1$ descents.  Replacing $k-i-1$ by $j$ yields
$$\lambda_n^{(
k)} = \sum_{j=0}^{k-1} {{n+k-j-1} \choose {k-j-1}}
(-1)^j \sum _{\sigma} sgn(\sigma) \sigma$$
which simplifies to $$\lambda_n^{(k)} = (-1)^{k-1} \sum_{j=0}^{k-1}
{{n+k-j-1} \choose n} \sum_{\sigma} sgn(\sigma) \sigma.$$In
each of
the last two equations, the sum is over $\sigma$ with $j$ descents.
Similarly,
$$\lambda_{n-1}^{(k)} = (-1)^{k-1} \sum_{j=0}^{k-1} {{n+k-j-2} \choose
{n-1}} \sum_{\sigma} sgn(\sigma) \sigma.$$
 
Lemma ~\ref{descent-increase} 
shows that for each $\sigma \in S_{n\setminus \{i\} }$ 
with $j$ descents, there are $j+1$ 
permutations $\tau \in S_n$ with $j$ descents such that $\partial
[i] \tau
= \pm \sigma$ and there are $n-j-1$ permutations $\tau \in S_n$ with 
$j+1$ descents such that $\partial [i]
\tau = \pm \sigma$.  When our  
boundary operator $\partial [i] $ deletes $\tau_d $ from $\tau $
to obtain $\sigma \in S_{n \setminus \{ i\} }$, we have 
$\partial [i] \tau = (-1)^{d-1} \sigma $, but Lemma ~\ref{sign-change}
proves that in this case, $sgn (\tau ) = (-1)^{i-d} sgn (\sigma )$.  
Combining these signs, observe that
$$\partial [i] sgn (\tau )\tau = (-1)^{i-1} sgn (\sigma )\sigma ,$$ 
independent of $d$.

Hence, the coefficient of $\sigma$ in $\partial [i] 
\lambda_n^{(k)}$ will be:
$$(-1)^{k-1}\cdot (-1)^{i-1} sgn(\sigma ) \left( {{n+k-j-1} \choose n} \cdot (j+1) + {{n+k-j-2} \choose n} \cdot (n-j-1) \right) $$which is equal to
$$(-1)^{k+ i - 2} sgn(\sigma) {{n+k-j-2} \choose n}\left( 
{{n+k-j-1} \over {k-j-1}} \cdot (j+1) + {{k-j-1} \over 
{k-j-1}}\cdot (n-j-1)  \right) . $$This simplifies to
$$(-1)^{k+i-2}sgn(\sigma) {{n+k-j-2} \choose n} \cdot {{kn} \over
{k-j-1}} = (-1)^{k+i-2} sgn(\sigma) {{n+k-j-2} \choose {n-1}} \cdot k.$$
This latter expression is the coefficient of $\sigma$ in 
$\lambda_{n-1}^{(k)}$ multiplied by 
$(-1)^{i-1}\cdot k$.  This holds for each $\sigma$, regardless of the
number of descents in $\sigma$, so we get 
$$\partial
 [i] \cdot \lambda_n^{(k)} = (-1)^{i-1}\cdot k \cdot 
\lambda_{n-1}^{(k)},$$ confirming identity (1).
 
We next prove that identity (1) implies identity (2).  Applying $\partial
 [i]$
to both sides of (2.3) gives:
$$(-1)^{k-1}(-1)^{i-1} \cdot k \cdot \lambda_{n-1}^{(k)} = 
\partial
[i]\sum_{j=1}^n k^j e_n^{(j)}.$$Applying (2.3) again, to the left hand
side yields,
$$(-1)^{i+k-2} \cdot k \cdot (-1)^{k-1} \sum_{j=1}^{n-1} k^j e_{n-1}^{(j)}
= \sum_{j=1}^n k^j \partial
 [i]e_n^{(j)}.$$
Hence,
$$(-1)^{i-1} \sum_{j=1}^{n-1} k^{j+1} e_{n-1}^{(j)} = \sum_{j=1}^n
k^j \partial
 [i] e_n^{(j)}.$$So,
$$0 = k \partial
 [i] e_n^{(1)} + \sum_{j=2}^n k^j 
\left(\partial
 [i] e_n^{(j)} - (-1)^{i-1} e_{n-1}^{(j-1)}\right) .$$
The fact that this holds for all $k$ implies that each coefficient of the
polynomial in $k$ is $0$, and so we get $\partial
 [i] e_n^{(1)} = 0$ and
$\partial
 [i] e_n^{(j)} = (-1)^{i-1} e_{n-1}^{(j-1)}$.  
\EOP

The second statement in Theorem 2.1 is particularly interesting when 
compared to a result that appears in the work of 
Gerstenhaber and Schack [GS].  In that work, the authors show that 
for the boundary $\delta$ in the usual complex for computing Hochschild 
homology of a commutative algebra, 
$$\delta e_n^{(k)} = e_{n-1}^{(k)} \delta$$
for all $n$ and $k$.  Note that this bears some similarity to the 
result we prove in Theorem 2.1 for the simplicial case although in the
simplicial case the boundary is applied on only one side and the 
Hodge index decreases by one rather than being constant.

For each $S \subseteq \underline{n}$, let $V_S$ denote the span of 
the $v_i$ for $i \in S$, and let $M_S$ denote the multilinear part of 
$V_S^{\otimes \mid S \mid}$.  We will continue using the Reiner-Webb
point of view so that the injective words on the set $S$ form a basis
for $M_S$.  
Note that $\partial(M_S) \subseteq \bigoplus_{i \in S} M_{S \backslash 
\{i\}}$, which means
we can decompose $\partial$ as a sum of the operators $\partial [i]$
for $i \in S$.  
 
Suppose $\mid S \mid = r$.  Then the symmetric group $S_r$ acts on $M_S$ 
by permuting the positions in which letters appear 
in the injective words on $S$.  For each $k$ with 
$1 \leq k \leq r$, let $M_S^{(k)}$ denote the image of $e_r^{(k)}$ under
this action.  By Theorem 2.1 (2), $\partial [i] (M_S^{(k)}) \subseteq
M_{S \backslash \{i\}
}$ for all $i \in S$.  So if we let $M_r^{(k)}$ denote
$\bigoplus_{ \mid S \mid = r} M_S^{(k)}$, then
$$\partial(M_r^{(k)}) \subseteq M_{r-1}^{(k-1)}.$$Thus, the complex
$(M_*,\partial_*)$ splits as a direct sum of the sub-complexes 
$C^{(k)}$ where $C^{(k)}$ is
$$(2.4)\ \ \ 0 \rightarrow M_n^{(k)} \rightarrow M_{n-1}^{(k-1)} 
\rightarrow \cdots \rightarrow M_{n-k+1}^{(1)} \rightarrow 0.$$Let
$H_*^{(k)}(M)$ denote the homology of the subcomplex $C^{(k)}$.  We
recall that $H_r^{(k)}(M)= 0$ unless $r=n$.  

Notice that the $S_n$-action on values which gives rise to the $S_n$-module
structure studied in [RW] commutes with the $S_r$ action on positions
in injective words in $M_r$.  Thus, it makes sense to study the 
$S_n$-module structure of $M_r^{(k)}$ for each $r$ and for
$H_n^{(k)}(M)$, with $S_n$ acting on values, despite the fact that the
Eulerian idempotents act on positions.  Our next result determines
$H_n^{(k)}(M)$ as an $S_n$-module.
To state this result, we will need some 
notation and results from [Ha].  
 
For each $\sigma \in S_n$, let $Z_{\sigma}$ denote the 
centralizer of $\sigma$.
In [Ha], a character $\chi_{\sigma}$ is defined as the induction of a linear
character $\Psi_{\sigma}$ from $Z_{\sigma}$ to $S_n$.  To describe $\Psi
_{\sigma}$, we first need a description of the $Z_{\sigma}$.  
Suppose $\sigma$
consists of $m_{\ell}$ $\ell$-cycles for each $\ell$.  Then $Z_{\sigma}$ 
is the direct product of $C_{\ell} wr S_{m_{\ell}}$ where $C_{\ell}$ is the
cyclic group of order $\ell$ and {\it wr} denotes wreath product.  
 
Let $\tau = \prod_{\ell} (\alpha_{\ell}; \beta_1, \beta_2, \ldots ,\beta_
{m_{\ell}})$ be an element of $Z_{\sigma}$ where $\alpha_{\ell} \in S_{m_
{\ell}}$ and each $\beta_i$ is in $C_{\ell}$.  Then
$$\Psi_{\sigma}(\tau) = \prod_{\ell} \prod_{i=1}^{m_{\ell}} \gamma_{\ell}
(\beta_i)$$where $\gamma_{\ell}$ is the linear character on $C_{\ell}$ which
assigns $e^{2 \pi i/\ell}$ to the generator of $C_{\ell}$.  
 
The following theorem from [Ha] will help us understand $M_r^{(k)}$.  

\begin{thm}[Hanlon]  \label{hanlon_hodge}
For each $n$ and $k$, let $I^{(k)}$ denote the
left ideal in {\bf C}$S_n$ generated by $e_n^{(k)}$.  As an $S_n$-module,
$$sgn * I_n
^{(k)} = \bigoplus_{\sigma} \chi_{\sigma}$$
where the sum is over a choice of representative from each conjugacy class 
that consists of permutations with exactly $k$ cycles.
\end{thm}

The $S_n$-modules $\bigoplus_{\sigma }\chi_{\sigma }$ in Theorem
~\ref{hanlon_hodge} have also appeared in a completely different context,
in work of Bergeron, Bergeron, Garsia on the free Lie algebra [BBG].  In 
contrast to Theorem ~\ref{hanlon_hodge} and [BBG], we will study 
$S_n$-modules in which we sum over conjugacy classes of derangements 
rather than conjugacy classes of permutations.

Our next result determines each $H_n^{(k)}(M)$ as an $S_n$-module, thereby 
providing a refinement of the theorem of Reiner and Webb which gives the
$S_n$-module structure of $H_n(M)$. 

\begin{thm}\label{hodge-piece}
For each $n$ and $k$,
$$sgn * H_n^{(k)}(M)  = \bigoplus_{\sigma} \chi_{\sigma}$$
where $*$ denotes the internal product, and the sum is over a 
choice of representative from each conjugacy class consisting of 
derangements with 
exactly $k$ cycles.  In particular, $\dim (H_n^{(k)}(M))$ equals the
number of derangements with exactly $k$ cycles.
\end{thm}

We will use cycle indices to prove Theorem ~\ref{hodge-piece}.  
For each $\sigma \in S_n$,
let $j_i(\sigma)$ denote the number of $i$-cycles of $\sigma$.  Let $a_1,
a_2, \ldots $ be a set of commuting indeterminates.  Define $Z(\sigma)$, the
cycle indicator of $\sigma$, to be 
$$Z(\sigma) = a_1^{j_1(\sigma)} a_2^{j_2 (\sigma)} a_3^{j_3(\sigma)} \cdots$$
So $Z(\sigma)$ is a monomial which identifies the cycle type of $\sigma$.
Thus, $\sigma$ and $\tau$ are conjugate in $S_n$ iff $Z(\sigma) = Z(\tau)$.
 
Let $\Psi$ be a class function on {\bf C}$S_n$.  The cycle index of $\Psi$
is 
$$Z(\Psi) = {1 \over {n!}} \sum_{\sigma \in S_n} \Psi(\sigma) Z(\sigma).$$
Since the monomial $Z(\sigma)$ uniquely identifies the conjugacy class of 
$\sigma$, two class functions are identical if and only if they have the 
same cycle index.
 
We will need two results about cycle indices from [Ha].  In the 
results below,
$\varepsilon_t$ denotes the trivial representation of $S_t$ (so 
$Z(\varepsilon_t) = {1 \over {t!}} \sum_{\sigma \in S_t} Z(\sigma)$) and
[] denotes the composition product on {\bf C}$[[a_1,a_2,\dots ]]^*$, i.e., 
for $A,B \in ${\bf C}$[[a_1,a_2,\ldots ]]$,
$$A[B] = A(a_i \leftarrow B(a_j \leftarrow a_{ij}))$$
where $\leftarrow$ denotes substitution.  Recall that $\mu $ denotes the 
ordinary number theory M\"obius function.
 
The following two results are proved in [Ha].

\begin{thm}[Hanlon]  Let $\sigma \in S_n$ with $Z(\sigma) = a_1^{j_1}
a_2 ^{j_2} \cdots a_n^{j_n}$.  Then
$$Z(\chi_{\sigma}) = \prod_{\ell = 1}^n Z(\varepsilon _{j_{\ell}}) 
\bigg[ {1 \over {\ell}} \sum_{d \mid \ell} \mu (d) a_d ^{\ell /d}\bigg] .$$
\end{thm}

\begin{thm}[Hanlon] Let $I_n^{(k)}$ denote the left ideal in {\bf C}$S_n$
generated by $e_n^{(k)}$.  Then
$$\sum _{n,k} Z(I_n^{(k)}) \lambda ^k = \prod_{\ell} 
(1 + (-1)^{\ell} a_{\ell})
^{{{-1}\over {\ell}} \sum_{d \mid \ell} \mu(d) \lambda^{\ell /d}}.$$
\end{thm}
We are now ready to give a proof of Theorem ~\ref{hodge-piece}.  

\proof
Recall that the Euler characteristic of a chain complex is the alternating
sum of the ranks of its homology groups, and the Hopf Trace Formula 
refines this to a statement about module structure.
Since the homology of $(M_*,\partial_r)$ 
vanishes except at the top degree, we deduce that
$$Z(H_n^{(k)}(M)) = \sum_{r=n-k+1}^n Z(M_r^{k-(n-r))})(-1)^{n-r}.$$
Note that 
$$M_r = \bigoplus_{\mid S \mid = r} M_S = ind_{S_r \times S_{n-r}}^{S_n} 
(Reg_r \otimes \varepsilon_{n-r})$$
where $Reg_r$ denotes the right-regular representation of $S_r$.  It follows
that 
$$M_r^{k-n+r} = ind_{S_r \times S_{n-r}}^{S_n} (I_r^{(k-n+r)} \otimes 
\varepsilon_{n-r}).$$
We will use one other well-known fact about cycle indices, namely that for
any virtual characters $\Psi$ of $S_r$ and $\Theta$ of $S_{n-r}$,
$$Z(ind_{S_r \times S_{n-r}}^{S_n}(\Psi \otimes \Theta)) = Z(\Psi)Z(\Theta).$$
Combining these facts we have:
\begin{eqnarray*}
\sum_n \sum_{k=1}^n Z(H_n^{(k)}(M))\lambda^k 
&=&\sum_n \sum_{k=1}^n \sum_{r=n-k+1}^n Z(M_r^{(k-(n-r))})(-1)^{n-r} 
\lambda^k \\
&=&\sum_n \sum_{k=1}^n \sum_{r=n-k+1}^n Z(I_r^{(k-n+r)})(\lambda^{k-n+r})
(-\lambda^{n-r})Z(\varepsilon_{n-r}) \\
&=&\bigg(\sum_{r,t} Z(I_r^{(t)})\lambda^t\bigg) \cdot \bigg(\sum_{m=0}^{\infty}
(-\lambda)^m Z(\varepsilon_m)\bigg) \\
&=&\bigg(\prod_{\ell} (1 + (-1)^{\ell} a_{\ell})^{{{-1}\over {\ell}} \sum_{d
\mid \ell} \mu(d) \lambda^{\ell / d}} \bigg) \cdot \exp\bigg(\sum_p 
{{(-\lambda)^p a_p} \over p}\bigg) 
\end{eqnarray*}
in the last step using the well-known fact that 
$$\sum_m Z(\varepsilon _m) = \exp\bigg(\sum_i {{a_i} \over i}\bigg) .$$
Thus,
$$\sum_n \sum_{k=1}^n Z(sgn *H_n^{(k)}(M))\lambda^k = F_1 \cdot F_2$$where
$$F_1 = \prod_{\ell} (1 - a_{\ell})^{{{-1} \over {\ell}} \sum_{d \mid \ell}
\mu(d) \lambda^{\ell/d}}$$and
$$F_2 = \exp \bigg( - \sum_p {{\lambda^p a_p} \over p}\bigg) .$$

We can rewrite $F_1$ as
\begin{eqnarray*}
F_1 &=& \prod_{\ell} \exp\bigg( \ln (1 - a_{\ell})\bigg( {{-1} \over {\ell}} 
\sum_{d \mid \ell} \mu(d) \lambda^{\ell /d}\bigg)\bigg) \\
&=& \prod_{\ell} \exp\bigg( \sum_{m=1}^{\infty} {{a_{\ell}^m} \over {m \ell}} 
\sum_{d \mid \ell} \mu(d) \lambda^{\ell / d}\bigg) . 
\end{eqnarray*}

Letting $p$ denote $\ell / d$ and $n$ denote $m d$, we have 

\begin{eqnarray*}
F_1 &=& \exp \bigg(\sum_{p,d,m} {1 \over {mpd}} \mu(d) \lambda^p 
  a_{dp}^m\bigg) \\
&=& \exp \bigg(\sum_p {{\lambda^p a_p} \over p}\bigg) 
  \bigg[\sum_{d,m} {1 \over {md}} \mu(d) a_d^m\bigg] \\
&=& \exp \bigg(\sum_p {{\lambda ^p a_p} \over p}\bigg)\bigg[
  \sum_n {1 \over n} \sum_{d \mid n} \mu(d) a_d^{n/d}\bigg] \\
&=& \exp \bigg(\sum_p {{\lambda^p a_p} \over p}\bigg) 
  \exp \bigg(\sum_p {{\lambda^p a_p} \over p}\bigg)
  \bigg[\sum_{\ell \ge 2} {1 \over {\ell}}
  \sum_{d \mid \ell} \mu(d) a_d^{\ell/d}\bigg] .
\end{eqnarray*}

So,
$$
\sum_{n,k} Z(sgn*H_n^{(k)}(M))\lambda^k = \exp 
  \bigg(\sum_p {{\lambda^p a_p} \over p}\bigg)
  \bigg[\sum_{\ell \ge 2} {1 \over {\ell}} \sum_{d \mid \ell} \mu(d) a_d^
  {\ell/d}\bigg] ,$$
which proves the result.
\EOP 

We conclude this section by showing how to recover Theorem 1.1 from 
Theorem 2.3.  Setting $\lambda = 1$ in Theorem 2.3, we obtain

\begin{eqnarray*}
\sum_n Z(sgn*H_n(M)) &=& \bigg(\sum_{n,k} Z(sgn*H_n^{(k)}(M) \bigg) \\
&=&\exp \bigg(\sum_p {{a_p} \over p}\bigg)\bigg[\sum_{\ell \ge 2} 
  {1 \over {\ell}} \sum_{d \mid {\ell}} \mu(d) a_d^{\ell /d} \bigg] \\
&=&\exp \bigg(\sum_p \sum_{\ell \ge 2} \sum_{d \mid \ell} {1 \over {p \ell}}
\mu(d) a_{dp}^{\ell /d} \bigg) 
\end{eqnarray*}

Letting $u = \ell/d$, we have
$$\sum_n Z(sgn *H_n(M)) = \exp\bigg(\sum_{p,d,u} {1 \over {pdu}} 
  \mu(d) a_{dp}^u\bigg) \cdot \exp\bigg( -\sum_p {{a_p} \over p}\bigg) $$
where the latter factor accounts for the provision that $\ell$ cannot 
equal $1$.  So, substituting $v$ for $pd$ yields

\begin{eqnarray*}
\sum_n Z(sgn*H_n(M)) &=& \exp \bigg(\sum_{u,v} {{a_v^u} \over {uv}} 
  (\sum_{d \mid v}\mu(d))\bigg) \cdot \exp \bigg( -\sum_p {{a_p} 
  \over p}\bigg)\\
&=&\exp \bigg(\sum_u {{a_1^u} \over u}\bigg) \cdot \exp\bigg(- 
  \sum_p {{a_p} \over p}\bigg) .
\end{eqnarray*}

Thus,
\begin{eqnarray*}
\sum_n Z(H_n(M)) &=& \bigg( {1 \over {1 - a_1}}\bigg)
  \cdot \exp\bigg(\sum_p {{(-1)^p a_p} \over p}\bigg) \\
&=& \bigg( \sum_k Z(Reg_k)\bigg) \cdot \bigg(\sum_{m=0}^{\infty} 
  (-1)^m Z(\varepsilon_m)\bigg) \\
&=& \sum_n \sum_k Z(ind_{S_{n-k} \times S_k}^{S_n} 
(\varepsilon_{n-k} \otimes Reg_k))(-1)^{n-k} 
\end{eqnarray*}

which is the Reiner-Webb Theorem.
 
\section{Signed random to random shuffles}

 In recent work, Uyemura-Reyes [Uy]
considers random to random shuffling of a deck of $n$ cards
and conjectures that the
transition matrix, when normalized to have integer entries, also has 
integer spectrum.  For small values of $n$, he notes that the nullspace
of the transition matrix has dimension equal to the number of derangements
of $n$ and that the nullspace, as an $S_n$-module, carries the same 
representation, up to a sign twist, as the representation that appears 
on the right-hand side of Theorem 1.1.  In this section, we explain this
phenomenon by studying the Laplacian $L$ of the complex $(M_*, \partial_*)$.
We show that if the normalized random-to-random shuffle operator has
integral spectrum (as conjectured in [Uy]), then the Laplacian on each 
chain group in the complex of injective words will also have integral
spectrum.

\begin{defn}
For each $r$, let $\upsilon_r$
and $\Upsilon_r$
be the 
elements of the group algebra {\bf C}$S_r$ given by:
$$\upsilon_r = r\cdot id + \sum_{u < v} (v,u,u+1,\ldots ,v-1)
+ \sum_{u>v} (v,u,u-1,\ldots,v+1)$$
and 
$$\Upsilon_r = r \cdot id + \sum_{u < v}(-1)^{v-u} (v,u,u+1, \ldots ,v-1)
+ \sum_{u>v} (-1)^{u-v} (v,u,u-1, \ldots ,v+1).$$
If we think of $S_r$ as acting on a deck of $r$ cards by permuting the
positions of the cards, then $\upsilon_r$ sums permutations which pick at
random two positions $v$ and $u$ and move the card in position $v$ to 
position $u$.  Thus $\upsilon_r$ is $r^2$ times the transition matrix
for random to random shuffling.  Note that $\Upsilon_r$ is simply 
$\upsilon_r$ twisted by the sign.  Thus we will refer to $\Upsilon_r$
as the
{\it signed random to random shuffle} element in {\bf C}$S_r$.  
\end{defn}

The following
conjecture appears in the dissertation of Uyemura-Reyes.

\begin{conj}[Uyemura-Reyes]
The eigenvalues of $\upsilon_n$ are (rational) integers.
\end{conj}
 
As in Section 2, we will use the collection $B_r$
of injective words of 
length $r$ on the alphabet $\{1,2,\ldots ,n\}$ as a basis for $M_r$.
Let $\delta_r :M_r \rightarrow M_{r+1}$ be the transpose of $\partial_{r+1}$
with respect to the inner products on $M_r$ and $M_{r+1}$ which have
$B_r$ and $B_{r+1}$ as orthonormal bases.  So if $D$ is the matrix for
$\partial_{r+1}$ with respect to the bases $B_r$ and $B_{r+1}$, then
$D^t$ is the matrix for $\delta_r$ with respect to the same bases.  Note
that $\delta_*$ is a coboundary on $M_*$.  We let $H^*(M)$ denote the 
cohomology with respect to this coboundary.
 
Let $\Lambda_r :M_r \rightarrow M_r$ be the Laplacian 
$$\Lambda_r = \delta_{r-1} \cdot \partial_r + \partial_{r+1} \cdot \delta_r.$$
We recall the well-known fact that a basis for the kernel of the 
Laplacian $\Lambda_r$ gives a simultaneous basis for $H_r(M)$ and 
$H^r(M)$.  

\begin{thm}
The Laplacian on the top-dimensional chain group satisfies 
$\Lambda_n = \Upsilon_n$.  
\end{thm}

\proof
To apply the coboundary $\delta_{n-1}$ to a basis element
$j_1 j_2 \ldots j_{n-1}$, we must sum over all sequences $i_1 i_2 
\ldots i_n$ with coefficient being the $\underline{j}, \underline{i}$
entry from $\partial_n$.  Since $\partial_n (i_1 \ldots i_n)$ is a sum
of terms of the form $\pm j_1 \ldots j_{n-1}$ where $j_1 \ldots j_{n-1}$
is obtained by deleting an entry from $i_1 \ldots i_n$, the $
\underline{j},\underline{i}$ entry of $\partial_n$ is $0$ unless 
$\underline{j}$ is a subsequence of $\underline{i}$.  It follows that
if $v$ is the single number in $\{1,2,\ldots,n\}$ which is missing from
$\{j_1,\ldots,j_{n-1}\}$ then
$$\delta_{n-1}(j_1 \ldots j_{n-1}) = (vj_1 \ldots j_{n-1}) 
-(j_1 v j_2 \ldots j_{n-1}) + (j_1 j_2 v j_3 \ldots j_{n-1}) - \cdots $$
So $\delta_{n-1} \partial_n$ is the operator which acts on a sequence 
$i_1 i_2 \ldots i_n$ by removing an element and re-inserting it in all
possible ways.  Moreover, if the removed element occupies position $u$
and it is re-inserted in position $v$ then the sign of that operation 
is $(-1)^{(u-1)+(v-1)} = (-1)^{v-u}$.  On the other hand, $\partial_{n+1}
\delta_n =0$.  It follows that $\Lambda_n$ is
equal to $\Upsilon_n$ which proves the result.
\EOP

\noindent
As noted above, Uyemura-Reyes conjectures that the spectrum of $\upsilon_n$
is integral from which it would follow that the spectrum of 
$\Lambda_n$ is integral.  We end this section by showing relating $\Lambda_r$
to $\Lambda_n$.  From this relationship one can deduce that if Conjecture 3.1
holds, then $\Lambda_r$ has integral spectrum for all $r$.

\begin{thm}
Let $i_1 \ldots i_r$ be a basis element of $M_r$.  Let
$A$ denote $\{ i_1, \ldots ,i_r\}$ and let $\overline{A}$ denote the 
complement of $A$ in $\{1,2,\ldots,n\}$.  Then,
$$\Lambda_r(i_1 \ldots i_r) = ((r+1)(n-r)I + \Upsilon_r + \sum_{a \in A,
b \in \overline{A}} 
(a,b))(i_1 \ldots i_r)$$where $\Upsilon_r$ is acting by permutation
of positions on $i_1 \ldots i_r$ whereas $(a,b)$ in the last summation is
acting by permuting the values of the $i_j$ within the set $\{1,2,\ldots,n\}$.
\end{thm}

\proof
We write $\Lambda_r(i_1 \ldots i_r)$ as a sum of three 
expressions
$$\Lambda_r (i_1 \ldots i_r) = X + Y + Z$$
where X is the sum of all terms in $\partial_{r+1} \delta_r (i_1 \ldots i_r)$
in which $\delta_r$ inserts an element $j$ of $\overline{A}$ in some position 
$u$ and then $\partial_{r+1}$ removes the same number $j$, where Y is the 
sum of all terms in $\delta_{r-1} \partial_r (i_1 \ldots i_r)$ in which
$\partial_r$ removes an element $j \in A$ and $\delta_{r-1}$ re-inserts
that same element $j$ and where Z is the remaining terms in $\Lambda_r
(i_1 \ldots i_r)$.  
 
Note that:
$$X = (r+1)(n-r)$$
and that:
$$Y = \Upsilon_r.$$
It will take some considerable effort now to analyze $Z$.  
 
The terms $\tau$ in $\partial_{r+1} \delta_r (i_1 \ldots i_r)$ that 
contribute to $Z$ are those in which an element $j$ from $\overline{A}$
is inserted into $i_1 \ldots i_r$ at some position $u$ by $\delta_r$
and then one of the $i_{\ell}$ is removed by $\partial_{r+1}$.  For
each such $\tau$, there is a corresponding term $\hat{\tau}$ in 
$\delta_{r-1} \partial_r (i_1 \ldots i_r)$ where $i_{\ell}$ is removed
first by $\partial_r$ and then $j$ is inserted in position corresponding 
to $u$ by $\delta_{r-1}$.  It is straightforward to check that $\tau = 
-\hat{\tau}$ and so these terms cancel.  
 
There is one circumstance in which this cancellation does not eliminate
every term.  These are the terms $\tau$ where $j$ is inserted immediately
behind $i_{\ell}$, i.e., where $u=\ell+1$.  In this case, the term 
$\hat{\tau}$ which should cancel $\tau$ is already committed to cancel the 
term $\tau '$ in which $j$ is inserted immediately in front of 
$i_{\ell}$.  
 
For $j \in \overline{A}$ and $i_{\ell} \in A$, the term in which $j$ is 
inserted immediately behind $i_{\ell}$ and then $i_{\ell}$ is deleted 
has sign $+1$ and is obtained by acting on $i_1 i_2 \ldots i_r$ with the
transposition $(j,i_{\ell}) \in S_n$.  The result follows.
\EOP

It is worth noting that the operator $\Upsilon_r$, acting on positions, 
commutes with the action of $S_n$ on words of length $r$ in $1,2, \ldots ,n$.
The first part of Theorem ~\ref{integral-spectrum} alternatively 
follows from [F] or from the shelling for $K_n$ in [BW].

\begin{thm}\label{integral-spectrum}
For $r<n$,
\begin{enumerate}
\item
$\Lambda_r$ is positive definite.
\item
If Conjecture 3.1 holds, then the spectrum of $\Lambda_r$ is 
integral. 
\end{enumerate}
\end{thm}

\proof
For this argument, it will be helpful to reconceptualize 
$M_r$.  Let $i_1 i_2 \ldots i_r$ be a basis element of $M_r$ and let
$\{j_1, \ldots ,j_{n-r}\} = \overline{A}$.  We will identify $i_1 \ldots i_r$
with 
$$[i_1 \ldots i_r] = {1 \over {(n-r)!}} \sum_{\sigma \in S_{n-r}} 
i_1 \ldots i_r j_{\sigma 1} j_{\sigma 2} \cdots j_{\sigma (n-r)} \in M_n.$$
The advantage this has is that the operator $\sum_{a \in A,b \in \overline{A}}
(a,b)$ whose action seemed to depend on the actual set $A$ can be redefined 
as the operator: 
$$\Gamma = \sum_{a \in \{1,\ldots,r\},b\in \{r+1,\ldots,n\}} (a,b)$$
where the permutation $(a,b)$ is acting now by permutation of positions.
So, $\Lambda_r = ((r+1)(n-r)I) + \Upsilon_r + \Gamma$.  
 
Let $\Omega = ((r+1)(n-r)I) + \Gamma$.  Note that $\Omega$ can be written 
as:
$$\Omega = ((r+1)(n-r)I) + T(1,n) - T(1,r) - T(r+1,n)$$
where $T(u,v) = \sum_{u \le a < b \le v} (a,b)$.  Recall that 
$$M_r = (Reg_r \otimes \varepsilon_{n-r}) \uparrow _{S_r \times S_{n-r}}
^{S_n}$$
where $Reg_r$ denotes the regular representation of $S_r$.  Therefore,
$$M_r = \bigoplus_{\alpha \vdash r} f_{\alpha} (S^{\alpha} \otimes 
\varepsilon_{n-r}) \uparrow _{S_r \times S_{n-r}}^{S_n}$$where $S^{\alpha}$
denotes the Specht module indexed by $\alpha$ and $f_{\alpha}$ is the number
of standard Young tableaux of shape $\alpha$.  
 
For
$x$ a square in row $i$ and column $j$ of a Ferrer's diagram of $\alpha$,
recall that $c_x$, the {\it content of x}, is 
$j-i$.  A well-known result from the 
representation theory of $S_n$ states that for a Specht module $S^{\lambda}$
with $\lambda \vdash n$, $T_n$ acts as the scalar $\sum_{x \in \lambda} c_x$.
It follows that for every $\lambda \vdash n$ which occurs in $(S^{\alpha} 
\otimes \varepsilon_{n-r}) \uparrow_{S_r \times S_{n-r}}^{S_n}$, the operator
$\Omega$ acts as the scalar:
$$(r+1)(n-r) + \sum_{x \in \lambda} c_x - \sum_{x \in \alpha} c_x 
- \sum_{x \in (n-r)} c_x$$
which simplifies to expression 3.2: 
$$(r+1)(n-r) + \sum_{x \in \lambda /\alpha} c_x - {{n-r} \choose 2}. 
\eqno(3.2)$$

We will make two observations based on this formula.  The first is that the
eigenvalues of $\Omega$ are integral.  Also, both $\Omega$ and $\Upsilon_r$
are easily seen to be diagonalizable.  Moreover, they commute.  It follows 
that the eigenvalues of $\Lambda_r = \Upsilon_r + \Omega$ can be written 
as sums of eigenvalues of $\Upsilon_r$ and $\Omega$.  However, $\Upsilon_r$
is conjugate to $\upsilon_r$ and hence has the same spectrum.
Thus, if Conjecture 3.1
holds, then all eigenvalues of $\Lambda_r$ are sums of
integers.  This proves part (b) of the theorem.
 
To prove that $\Lambda_r$ is positive definite, first note that it is enough
to show that $\Omega$ is positive definite since $\Upsilon_r$ is positive 
semi-definite, being a direct sum ${n \choose r}$ copies of the Laplacian 
in top degree for the case with $n=r$.  To see that $\Omega$ is positive
definite, we start with the expression for the action of $\Omega$ on copies
of $S_n$ irreducibles given in (3.2) above.  The first observation follows
from the fact that $S^{\lambda}$ has nonzero multiplicity in 
$(S^{\alpha} \otimes \varepsilon_{n-r}) \uparrow _{S_r \times S_{n-r}}^{S_n}$
if and only if $\lambda / \alpha$ is a horizontal strip.  
Let $(\rho_1,\gamma_1), (\rho_2,\gamma_2), \ldots ,(\rho_{n-r},\gamma_{n-r})$
be the coordinates of the squares in $\lambda / \alpha$.  The fact that 
$\lambda / \alpha$ is a horizontal strip implies that 
$1 \le \gamma_1 < \gamma_2 \cdots < \gamma_{n-r}$.  Thus, 
$$\sum_s \gamma_s
\ge {{n-r} \choose 2} + (n-r).$$
Also, observe that if a square x of the Ferrer's diagram of $\lambda / \alpha$
is in row $i$, then there are $(i-1)$ squares of $\alpha$ in the rows above 
it.  So, $\sum (\rho_s-1) \le r$, i.e., $\sum \rho_s \le r + (n-r)$.  
 
Putting these bounds together gives that the eigenvalue $\omega$ given in 
formula (3.2) satisfies
$$\omega \ge (r+1)(n-r) + {{n-r} \choose 2} + (n-r) - (r + (n-r)) - 
{{n-r} \choose 2}$$
which simplifies to:
$$\omega \ge (r+1)(n-r) -r > 0.$$
\EOP

\end{document}